\documentstyle[12pt]{article}
\def\mineappendix{
        \setcounter{section}{1}
        \setcounter{subsection}{0}
        \def\thesection{\Alph{section}}
        \def\sectionap{\@startsection  {section}{1}{\z@}
                        {-3.5ex plus-1ex minus-.2ex} {0ex plus.2ex}
                        {\reset@font\Large\bf  Appendix:  \, }
                        }
        }
\makeatother
\def\Proclaim #1. #2\par{\bigbreak\noindent{\sc#1.\enspace}{\it#2}\par}

\font\Bbbfont=msbm10
\newfam\msbfam
\textfont\msbfam=\Bbbfont
\def\Bbb#1{{\fam\msbfam\relax#1}}
\newcommand{\eqref}[1]{equation~(\ref{#1})}
\newcommand{\Eqref}[1]{Equation~(\ref{#1})}

\newcommand{\gwii}[1]{\left< \hspace{-2pt} \left< \, #1 \,
        \right>  \hspace{-2pt} \right>_{0}}

\newcommand{\gwiione}[1]{\left< \hspace{-2pt} \left< \, #1 \,
        \right> \hspace{-2pt} \right>_{1}}

\newcommand{\gwiitwo}[1]{\left< \hspace{-2pt} \left< \, #1 \,
        \right> \hspace{-2pt} \right>_{2}}

\newcommand{\gwig}[1]{\left< \, #1 \, \right>_{g}}
\newcommand{\gwiig}[1]{\left< \hspace{-2pt} \left< \, #1 \,
    \right> \hspace{-2pt} \right>_{g}}

\newcommand{\grav}[2]{\tau_{#1}(\gamma_{#2})}
\newcommand{\grava}[1]{\tau_{#1}(\gamma_{\alpha})}

\newcommand{\ga}{\gamma_{\alpha}}
\newcommand{\gua}{\gamma^{\alpha}}
\newcommand{\gb}{\gamma_{\beta}}
\newcommand{\gub}{\gamma^{\beta}}
\newcommand{\gm}{\gamma_{\mu}}
\newcommand{\gum}{\gamma^{\mu}}

\newcommand{\vs}{{\cal S}}

\newcommand{\vd}{{\cal D}}
\newcommand{\vw}{{\cal W}}
\newcommand{\vv}{{\cal V}}

\newcommand{\qp}{\bullet}

\newcommand{\qpcc}[2]{{\cal #1} \bullet {\cal #2}}

\newtheorem{lem}{Lemma}[section]

\newtheorem{thm}[lem]{Theorem}

\title{Relations Among Universal Equations For Gromov-Witten Invariants}
\author{Xiaobo Liu \thanks{Research partially supported by
            Alfred P. Sloan Research Fellowship and NSF research grant}}
\date{}

\begin{document}
\maketitle

It is well known that relations in the tautological ring of moduli
spaces of pointed stable curves give partial differential
equations for Gromov-Witten invariants of compact symplectic
manifolds. These equations do not depend on the target symplectic
manifolds and therefore are called universal equations for
Gromov-Witten invariants. In the case that the quantum cohomology
of the symplectic manifolds are semisimple, it is expected that
higher genus Gromov-Witten invariants are completely determined by
such universal equations and genus-0 Gromov-Witten invariants.
This has been proved for genus-1 (cf. \cite{DZ}) and genus-2 (cf.
\cite{L2}) cases.
Universal equations also play very important role in the understanding of
the Virasoro conjecture (cf. \cite{EHX}). The genus-0 Virasoro conjecture
for all compact symplectic manifolds follows from a universal equation called
the genus-0 topological recursion relation (cf. \cite{LT}).
For projective varieties, we expect that such universal equations
reduce higher genus Virasoro conjecture to an $SL(2)$
symmetry for the generating function of the Gromov-Witten invariants.
Again this has been proved for genus-1 (\cite{L1}) and genus-2 (\cite{L2})
cases.

In this paper, we will discuss the relation among known universal equations
for Gromov-Witten invariants. We hope that the understanding of such relations
would be helpful to the study of both Gromov-Witten invariants and the topology
of the moduli spaces of pointed curves. Relations among genus-2 universal equations
were studied in \cite{L2}. It was proved that the three universal equations
in \cite{G2} and \cite{BP} implies certain complicated genus-1 relations
(see equations (\ref{eqn:A1A2}) and (\ref{eqn:A1A2B}) below). Modulo these genus-1 relations,
the three known genus-2 equations can be reduced to only two equations. The main result
of this paper is that the genus-1 relations derived in \cite{L2} follow
from a known genus-1 relation found in \cite{G1} and the genus-0 
and genus-1 topological recursion relations.
This completes the discussion of relations among genus-2 equations in \cite{L2}.

Part of the work in this paper was done while the author visited Max-Planck Institute
for Mathematics at Bonn during the workshop on Frobenius Manifolds, Quantum Cohomology,
and Singularity Theory. The author would like to thank the organizers of the workshop
for invitation. He would also like to thank E. Getzler for conversations related to this
work.

\section{Generating functions for Gromov-Witten invariants}

Let $M$ be a compact symplectic manifolds. For simplicity, we assume
$H^{\rm odd}(M; {\Bbb C}) = 0$.
The {\it big phase} is by definition the product of infinite copies of
$H^{*}(M; {\Bbb C})$, i.e.
\[  P := \prod_{n=0}^{\infty} H^{*}(M; {\Bbb C}). \]
Fix a basis $\{ \gamma_{1}, \ldots, \gamma_{N} \}$ of
$H^{*}(M; {\Bbb C})$ with $\gamma_{1} = 1$ be the identity of the ordinary
cohomology ring of $M$. Then we denote the corresponding basis for
the $n$-th copy of $H^{*}(M; {\Bbb C})$ in $P$ by
$\{\tau_{n}(\gamma_{1}), \ldots, \tau_{n}(\gamma_{N}) \}$.
We call $\grava{n}$ a {\it descendant} of $\gamma_{\alpha}$ with descendant
level $n$.
We can think of $P$ as an infinite dimensional vector space with basis
$\{ \grava{n} \mid 1 \leq \alpha \leq N, \, \, \, n \in {\Bbb Z}_{\geq 0} \}$
where ${\Bbb Z}_{\geq 0} = \{ n \in {\Bbb Z} \mid n \geq 0\}$.
Let
$(t_{n}^{\alpha} \mid 1 \leq \alpha \leq N, \, \, \, n \in {\Bbb Z}_{\geq 0})$
be the corresponding coordinate system on $P$.
For convenience, we identify $\grava{n}$ with the coordinate vector field
$\frac{\partial}{\partial t_{n}^{\alpha}}$ on $P$ for $n \geq 0$.
If $n<0$, $\grava{n}$ is understood as the $0$ vector field.
We also abbreviate $\grava{0}$ as $\gamma_{\alpha}$.
Any vector field of the form $\sum_{\alpha} f_{\alpha} \ga$, where $f_{\alpha}$
are functions on the big phase space, is called a {\it primary vector field}.
We use $\tau_{+}$ and $\tau_{-}$ to denote the operator which shift the level
of descendants, i.e.
\[ \tau_{\pm} \left(\sum_{n, \alpha} f_{n, \alpha} \grava{n}\right)
    = \sum_{n, \alpha} f_{n, \alpha} \grava{n \pm 1} \]
where $f_{n, \alpha}$ are functions on the big phase space.

We will use the following {\it notational conventions}:
Lower case Greek letters, e.g. $\alpha$, $\beta$, $\mu$, $\nu$,
$\sigma$,..., etc., will be used to index the cohomology classes.
The range of these indices is from $1$ to $N$.
Lower case English
letters, e.g. $i$, $j$, $k$, $m$, $n$, ..., etc., will be used to
index the level of descendants. Their range is the set of all
non-negative integers, i.e. ${\Bbb Z}_{\geq 0}$. All summations are
over the entire ranges of the indices unless otherwise indicated.
Let
\[ \eta_{\alpha \beta} = \int_{V} \gamma_{\alpha} \cup
    \gamma_{\beta}
\]
 be the intersection form on $H^{*}(V, {\Bbb C})$.
We will use $\eta = (\eta_{\alpha \beta})$ and $\eta^{-1} =
(\eta^{\alpha \beta})$ to lower and raise indices.
For example,
\[ \gua :=  \eta^{\alpha \beta} \gb.\]
Here we are using the summation convention that repeated
indices (in this formula, $\beta$) should be summed
over their entire ranges.

Let
\[ \gwig{\grav{n_{1}}{\alpha_{1}} \, \grav{n_{2}}{\alpha_{2}} \,
    \ldots \, \grav{n_{k}}{\alpha_{k}}} \]
be the genus-$g$ descendant Gromov-Witten invariant associated
to $\gamma_{\alpha_{1}}, \ldots, \gamma_{\alpha_{k}}$ and nonnegative
integers $n_{1}, \ldots, n_{k}$
(cf. \cite{W}, \cite{RT}, \cite{LiT}, and \cite{M}).
Let $\overline{\cal M}_{g, k}(M; d)$ be the moduli space of stable maps
from genus-g k-pointed curves to $M$ with degree $d \in H_{2}(M; {\Bbb Z})$,
and $\Psi_{i}$ the first Chern class of the tautological line bundle
over $\overline{\cal M}_{g, k}(M; d)$ whose geometric fiber over a stable map
is the cotangent space of the domain curve at $i$-th marked point.
Then the integers $n_{1}, \ldots, n_{k}$ above represent powers of
$\Psi_{1}, \ldots, \Psi_{k}$.
The genus-$g$
generating function is defined to be
\[ F_{g} =  \sum_{k \geq 0} \frac{1}{k!}
         \sum_{ \begin{array}{c}
        {\scriptstyle \alpha_{1}, \ldots, \alpha_{k}} \\
                {\scriptstyle  n_{1}, \ldots, n_{k}}
                \end{array}}
                t^{\alpha_{1}}_{n_{1}} \cdots t^{\alpha_{k}}_{n_{k}}
    \gwig{\grav{n_{1}}{\alpha_{1}} \, \grav{n_{2}}{\alpha_{2}} \,
        \ldots \, \grav{n_{k}}{\alpha_{k}}}. \]
This function is understood as a formal power series of
$t_{n}^{\alpha}$.

Introduce
a $k$-tensor
 $\left< \left< \right. \right. \underbrace{\cdot \cdots \cdot}_{k}
        \left. \left. \right> \right> $
defined by
\[ \gwiig{{\cal W}_{1} {\cal W}_{2} \cdots {\cal W}_{k}} \, \,
         := \sum_{m_{1}, \alpha_{1}, \ldots, m_{k}, \alpha_{k}}
                f^{1}_{m_{1}, \alpha_{1}} \cdots f^{k}_{m_{k}, \alpha_{k}}
        \, \, \, \frac{\partial^{k}}{\partial t^{\alpha_{1}}_{m_{1}}
            \partial t^{\alpha_{2}}_{m_{k}} \cdots
            \partial t^{\alpha_{k}}_{m_{k}}} F_{g},
 \]
for (formal) vector fields ${\cal W}_{i} = \sum_{m, \alpha}
        f^{i}_{m, \alpha} \, \frac{\partial}{\partial t_{m}^{\alpha}}$ where
$f^{i}_{m, \alpha}$ are (formal) functions on the big phase space.
We can also view this tensor as the $k$-th covariant derivative
of $F_{g}$. This tensor is called the {\it $k$-point
(correlation) function}.

For any vector fields $\vw_{1}$ and $\vw_{2}$ on the big phase space,
the quantum product of $\vw_{1}$ and $\vw_{2}$ is defined by (cf. \cite{L2})
\[ \vw_{1} \qp \vw_{2} := \gwii{\vw_{1} \, \vw_{2} \, \gua} \ga. \]
Define
\[ T(\vw) := \tau_{+}(\vw) - \gwii{\vw \, \gua} \ga \]
for any vector field $\vw$. The operator $T$
was introduced in \cite{L2} as a convenient tool in the study of topological recursion
relations. Let $\overline{\cal M}_{g, k}$ be the moduli space of genus-g stable
curves with $k$-marked points, and $\psi_{i}$ the first Chern class of the
tautological line bundle over $\overline{\cal M}_{g, k}$ whose geometric fiber
over a stable curve is the cotangent space of the curve at the $i$-th marked point.
 When we translate a relation in the tautological
ring of $\overline{\cal M}_{g, k}$ to differential equations for
generating functions of Gromov-Witten invariants, the $\psi$ classes
correspond to the operator $T$. The reason for such correspondence lies in
the relation between the classes
$\Psi_{i} \in H^{2}(\overline{\cal M}_{g, k}(M, d); \, {\Bbb Q})$
and $\psi_{i} \in H^{2}(\overline{\cal M}_{g, k}; \, {\Bbb Q})$.
There is a canonical map $St: \overline{\cal M}_{g, k}(M, d) \longrightarrow
\overline{\cal M}_{g, k}$ which forgets the map to the target manifold $M$ and stabilizes
the domain curve for each element of $\overline{\cal M}_{g, k}(M, d)$. The difference
$\Psi_{i} - St^{*} (\psi_{i})$ is represented by a cycle containing elements whose domain curves
consist of one genus-$g$ and one genus-0 components, with the $i$-th marked point lying
on the genus-0 component (cf. \cite{G2} and \cite{KM}).
When we interpret this relation for Gromov-Witten invariants,
$\Psi_{i}$ corresponds to $\tau_{+}$, the $i$-th marked point is associated with a cohomology
class, which will be extended to a vector field $\vw$ by linearity, and
the node joining the two irreducible components are associated with the diagonal
class $\gua \otimes \ga$ in $H^{*}(M \times M; {\Bbb C})$. Consequently, $St^{*}(\psi_{i})$
must correspond to $T$.

The operator $T$ also has an algebraic interpretation.
Let \[\vs := - \sum_{m, \alpha} \tilde{t}^{\alpha}_{m}
        \grav{m-1}{\alpha} \]
be the {\it string vector field}, where
\[ \tilde{t}^{\alpha}_{m}:= t_{m}^{\alpha}
        - \delta_{m, 1} \delta_{\alpha, 1}. \]
By the second derivative of the string equation (see Section~\ref{sec:StringDilaton}),
the operator $T$ can be rewritten as
\[ T(\vw) = \tau_{+}(\vw) - \vs \qp \tau_{+}(\vw). \]
Therefore $T$ measures the difference between $\vs$ and the identity of the quantum product
on the big phase space, which actually does not exist by our definition of the quantum
product.

Let $\nabla$ be the trivial flat connection on the big phase space with respect
to which $\grava{n}$ are parallel vector fields for all $\alpha$ and $n$.
Then the covariant derivative of the quantum product is given by
\begin{equation} \label{eqn:derQP}
\nabla_{\vw_{3}} (\vw_{1} \qp \vw_{2})
    = (\nabla_{\vw_{3}} \vw_{1}) \qp \vw_{2}
        +  \vw_{1} \qp (\nabla_{\vw_{3}} \vw_{2})
        + \gwii{\vw_{1} \, \vw_{2} \, \vw_{3} \, \gua} \ga
\end{equation}
and the covariant derivative of the operator $T$ is given by
\begin{equation} \label{eqn:derT}
 \nabla_{\vw_{2}} \, \, T(\vw_{1}) = T(\nabla_{\vw_{2}} \vw_{1})
    - \vw_{2} \qp \vw_{1}
\end{equation}
for any vector fields $\vw_{1}, \vw_{2}$ and  $\vw_{3}$
(cf. \cite[Equation (8) and Lemma 1.5]{L2}).

\section{Universal equations}
\label{sec:univ}

We will write universal equations of Gromov-Witten invariants as
equations among tensors on the big phase space defined by
generating functions of Gromov-Witten invariants.
The set of tensors whose vanishing are equivalent to some universal equations
form an ideal of
contra-variant tensor algebra on the big phase space. We call this ideal the
{\it ideal of universal relations} and denote
it by $U$. This ideal is closed under the covariant differentiation
defined by $\nabla$, and therefore is a differential ideal.

If $T_{1}, \ldots, T_{k}$ are $k$ contra-variant tensors, we use
$I_{a}[T_{1}, \cdots, T_{k}]$ to denote the
ideal algebraically generated by $T_{1}, \ldots, T_{k}$, i.e.
$I_{a}[T_{1}, \cdots, T_{k}]$ is the set of tensors of the form
$\sum_{i=1}^{k} S^{'}_{i} \otimes T_{i} + T_{i} \otimes S^{''}_{i}$
where $S^{'}_{i}$ and $S_{i}^{''}$ are arbitrary contra-variant tensors
on the big phase space. We also use
$I_{d}[T_{1}, \cdots, T_{k}]$ to denote the
ideal differentially generated by $T_{1}, \ldots, T_{k}$, i.e.
the smallest ideal which contains $T_{1}, \ldots, T_{k}$ and is closed
under the covariant differentiation.

In the theory of Gromov-Witten invariants, we often need to prove two complicated
expressions are equal by using certain universal equations. Since universal equations
of Gromov-Witten invariants are very complicated, it is very cumbersome to write out
explicitly derivations in each step. For convenience of presentation, we introduce the
following phrases: If  one expression can be derived from another one by
using certain universal equations $T_{1} = 0, \cdots, T_{k}=0$, but not using
their covariant derivatives, then we say that these {\it two expressions are equal
 modulo $I_{a}[T_{1}, \cdots, T_{k}]$}. If in the derivation, we also need covariant
 derivatives of the corresponding universal equations, we say that these {\it two expressions
 are equal modulo $I_{d}[T_{1}, \cdots, T_{k}]$}.

\subsection{genus-0 universal relations}

\label{sec:g0univ}

Define
\[ \rho_{0}(\vw_{1}, \vw_{2}, \vw_{3})
        := \gwii{T(\vw_{1}) \, \vw_{2} \, \vw_{3}} \]
where $\vw_{i}$  are vector fields on the big phase space.
The genus-0 topological recursion relation is equivalent to
$\rho_{0} = 0$. Hence $\rho_{0} \in U$, and so does the tensor $C_{0}$
defined by
\[ C_{0}(\vw_{1}, \vw_{2}, \vw_{3}, \vw_{4}) :=
    (\nabla_{\vw_{3}} \rho_{0})(\vw_{1}, \vw_{2}, \vw_{4})
    -(\nabla_{\vw_{2}} \rho_{0})(\vw_{1}, \vw_{3}, \vw_{4}). \]
 Using \eqref{eqn:derT}, it is straightforward to show
that
\[ C_{0}(\vw_{1}, \vw_{2}, \vw_{3}, \vw_{4}) =
        \gwii{ (\vw_{1} \qp \vw_{2}) \, \vw_{3} \, \vw_{4}}
        - \gwii{\vw_{1} \, (\vw_{2} \qp \vw_{3}) \, \vw_{4}} \]
for any vector fields $\vw_{1}, \ldots, \vw_{4}$.
Therefore $C_{0} = 0$ implies that the quantum product "$\qp$" on the
big phase space is associative.

\subsection{genus-1 universal relations}
\label{sec:g1univ}

The genus-1 topological recursion relation is equivalent to
$\rho_{1} = 0$ where $\rho_{1}$ is defined by
\[ \rho_{1}(\vw) = \gwiione{T(\vw)}
        - \frac{1}{24} \gwii{\vw \, \gua \, \ga}. \]
Define
\begin{eqnarray*}
 G(\vw_{1}, \vw_{2}, \vw_{3}, \vw_{4})
& := &
        \sum_{g \in S_{4}} \left\{
                3 \gwiione{ \{\vw_{g(1)} \qp \vw_{g(2)} \}
                        \, \{\vw_{g(3)} \qp \vw_{g(4)} \} }
              \right. \\
    && \hspace{25pt}
        -  4 \gwiione{ \{\vw_{g(1)} \qp \vw_{g(2)} \qp
                        \vw_{g(3)} \} \, \vw_{g(4)} }  \\
        && \hspace{25pt} -  \gwii{ \{ \vw_{g(1)} \qp \vw_{g(2)} \} \,
                                \vw_{g(3)} \, \vw_{g(4)} \,
                                        \gamma^{\alpha} }
        \gwiione{ \gamma_{\alpha} }   \\
    && \hspace{25pt}
         +  2 \gwii{ \vw_{g(1)} \, \vw_{g(2)} \, \vw_{g(3)} \,
                                        \gamma^{\alpha} }
        \gwiione{ \{\gamma_{\alpha} \bullet \vw_{g(4)} \} } \\
    && \hspace{25pt}
        + \frac{1}{6} \gwii{ \vw_{g(1)}
                \, \vw_{g(2)} \, \vw_{g(3)} \,
                                        \gamma^{\alpha} }
                \gwii{ \gamma_{\alpha} \, \vw_{g(4)} \, \gamma_{\beta}
                                     \,   \gamma^{\beta} }
                \\
        && \hspace{25pt}
        + \frac{1}{24} \gwii{ \vw_{g(1)} \,
            \vw_{g(2)} \, \vw_{g(3)} \, v_{g(4)} \,
                                        \gamma^{\alpha} }
        \gwii{ \gamma_{\alpha} \, \gamma_{\beta} \,
                                        \gamma^{\beta} }
                        \\
        && \hspace{25pt} \left.
        - \frac{1}{4} \gwii{ \vw_{g(1)} \,  \vw_{g(2)} \,
                        \gamma^{\alpha} \, \gamma^{\beta} }
                \gwii{ \gamma_{\alpha} \, \gamma_{\beta} \,
                                \vw_{g(3)} \, \vw_{g(4)} }
                        \right\}
\end{eqnarray*}
for any vector fields $\vw_{1}, \ldots, \vw_{4}$. Then the genus-1 equation
discovered by Getzler \cite{G1} is equivalent to $G=0$.
So $\rho_{1} \in U$ and $G \in U$.

\subsection{genus-2 universal relations}

Define
\begin{eqnarray*}
&& \rho_{2; 1} (\vw) := \gwiitwo{T^{2}(\vw)} - A_{1}(\vw), \hspace{200pt}  \\
&&
\rho_{2; 2} (\vw, \vv) :=
    \gwiitwo{T({\cal W}) \, T({\cal V})} - 3 \gwiitwo{ T(\qpcc{W}{V})}
    - A_{2}(\vw, \vv), \\
&& \rho_{2; 3}(\vw_{1}, \vw_{2}, \vw_{3}) :=
 2 \gwiitwo{\{{\cal W}_{1} \bullet {\cal W}_{2} \bullet {\cal W}_{3} \}}
- 2 \gwii{{\cal W}_{1} \, {\cal W}_{2} \, {\cal W}_{3} \, \gua}
    \gwiitwo{T(\ga)}
\nonumber \\
&& \hspace{60pt} + \sum_{g \in S_{3}}
    \gwiitwo{{\cal W}_{g(1)} \,
        T({\cal W}_{g(2)} \qp {\cal W}_{g(3)})}
    -\gwiitwo{ T({\cal W}_{g(1)}) \, \{{\cal W}_{g(2)}
        \qp {\cal W}_{g(3)}\}}
\nonumber \\
&& \hspace{60pt} \, \, - \, \, B(\vw_{1}, \vw_{2}, \vw_{3}),
\end{eqnarray*}
where
\begin{eqnarray*}
A_{1}(\vw) &=&
    \frac{7}{10} \gwiione{\ga} \gwiione{\{\gua \bullet \vw\}}
    + \frac{1}{10} \gwiione{\ga \, \{\gua \bullet \vw\}} \\
&&     - \frac{1}{240} \gwiione{\vw \, \{\ga \bullet \gua\}}
    + \frac{13}{240} \gwii{\vw \, \ga \, \gua \, \gub}
        \gwiione{\gb} \\
&&    + \frac{1}{960} \gwii{ \vw \, \gua \, \ga \, \gub \, \gb}, \\
A_{2}(\vw, \vv) &=&
     \frac{13}{10} \gwii{\vw \, \vv \, \gua \, \gub}
        \gwiione{\ga} \gwiione{\gb}
    + \frac{4}{5} \gwiione{\vw \, \gua}
            \gwiione{\{\ga \bullet \vv\}} \\
&&    + \frac{4}{5}
        \gwiione{\vv \, \gua} \gwiione{\{\ga \bullet \vw\}}
    - \frac{4}{5}
        \gwiione{ \{\vw \bullet \vv\}  \, \gua} \gwiione{\ga} \\
&&    + \frac{23}{240} \gwii{ \vw \, \vv \,
        \gua \, \ga \, \gub} \gwiione{\gb}
    + \frac{1}{48}\gwii{ \vw \, \gua \, \ga \, \gub}
        \gwiione{ \gb \,\vv}  \\
&&    + \frac{1}{48}\gwii{ \vv \, \gua \, \ga \, \gub}
        \gwiione{ \gb \,\vw}
    - \frac{1}{80} \gwiione{\vw \, \vv \, \{\gua \bullet \ga\}} \\
&&    + \frac{7}{30} \gwii{\vw \, \vv \, \gua \, \gub}
        \gwiione{ \ga \, \gb}
    + \frac{1}{30}
        \gwiione{\ga \, \{\gua \bullet \vw\} \, \vv} \\
&&    + \frac{1}{30}
        \gwiione{\ga \, \{\gua \bullet \vv\} \, \vw}
    - \frac{1}{30}
        \gwiione{\{ \vw \bullet \vv\} \ga \, \gua} \\
&&    + \frac{1}{576}\gwii{\vw \, \vv \, \gua \, \ga
        \, \gub \, \gb},
\end{eqnarray*}
and
\begin{eqnarray*}
&& B(\vw_{1}, \vw_{2}, \vw_{3}) \\
&=& \frac{1}{5} \gwii{\vw_{1} \, \vw_{2} \, \vw_{3} \, \, \gua \, \gub}
    \gwiione{\ga} \gwiione{\gb}
   - \frac{6}{5} \gwii{\vw_{1} \, \vw_{2} \, \vw_{3} \, \gua}
    \gwiione{\ga \, \gub} \gwiione{\gb}
    \nonumber \\
&&
+ \frac{1}{120} \gwii{\vw_{1} \,   \vw_{2} \, \vw_{3}
         \, \gua \, \ga \, \gub}
    \gwiione{\gb}
- \frac{1}{120} \gwiione{\vw_{1} \, \vw_{2} \, \vw_{3}
       \, \{ \gua \bullet \ga \}}
    \nonumber \\
&&
+ \frac{1}{10} \gwii{\vw_{1} \, \vw_{2} \, \vw_{3}
         \, \gua \, \gub}
    \gwiione{\ga \, \gb}
- \frac{1}{20} \gwii{\vw_{1} \, \vw_{2} \, \vw_{3} \, \gua}
    \gwiione{\ga \, \gub \, \gb}
    \nonumber \\
&&
- \frac{1}{5} \sum_{\sigma \in S_{3}}
     \gwii{\vw_{\sigma(1)} \, \vw_{\sigma(2)} \, \gua \, \gub}
   \gwiione{\ga} \gwiione{\gb \, \, \vw_{\sigma(3)}}
    \nonumber \\
&&
+ \frac{2}{5}\sum_{\sigma \in S_{3}}
    \gwiione{\{\vw_{\sigma(1)} \bullet \ga \}}
    \gwiione{\gua \, \, \vw_{\sigma(2)} \,
        \, \vw_{\sigma(3)}}
    \nonumber \\
&&
- \frac{3}{5} \sum_{\sigma \in S_{3}}
    \gwiione{\vw_{\sigma(1)} \, \, \{ \vw_{\sigma(2)} \bullet \gua \}}
    \gwiione{\ga \, \, \vw_{\sigma(3)}}
    \nonumber \\
&&
+ \frac{3}{10}\sum_{\sigma \in S_{3}}
    \gwiione{\vw_{\sigma(1)} \, \, \ga}
    \gwiione{\gua \, \{\vw_{\sigma(2)} \bullet \vw_{\sigma(3)} \}}
    \nonumber \\
&&
- \frac{1}{5}\sum_{\sigma \in S_{3}}
    \gwiione{\ga}
    \gwiione{\gua \, \, \vw_{\sigma(1)} \, \,
        \{\vw_{\sigma(2)} \bullet \vw_{\sigma(3)}\}}
    \nonumber \\
&&
- \frac{1}{80} \sum_{\sigma \in S_{3}}
    \gwii{\vw_{\sigma(1)} \, \,
        \vw_{\sigma(2)} \, \, \gua \, \ga \, \gub}
    \gwiione{\gb \, \, \vw_{\sigma(3)}}
    \nonumber \\
&&
+ \frac{1}{80} \sum_{\sigma \in S_{3}}
    \gwii{\vw_{\sigma(1)} \, \,
            \gua \, \ga \, \gub}
    \gwiione{\gb \, \, \vw_{\sigma(2)} \, \,
            \vw_{\sigma(3)}}
    \nonumber \\
&&
- \frac{1}{20} \sum_{\sigma \in S_{3}}
     \gwiione{\vw_{\sigma(1)} \, \, \ga \, \gb}
        \gwii{\gua \, \gub \, \, \vw_{\sigma(2)} \,
        \, \vw_{\sigma(3)}}
    \nonumber \\
&&
+ \frac{1}{60} \sum_{\sigma \in S_{3}}
    \gwiione{\vw_{\sigma(1)} \, \,
            \vw_{\sigma(2)} \, \, \ga \, \{\gua \bullet \vw_{\sigma(3)} \}}
    \nonumber \\
&&
- \frac{1}{120} \sum_{\sigma \in S_{3}}
     \gwiione{\vw_{\sigma(1)} \, \,
        \gua \, \ga \, \{ \vw_{\sigma(2)} \bullet \vw_{\sigma(3)} \} }.
\end{eqnarray*}
Note that the 3 complicated tensors $A_{1}, A_{2}, B$ are symmetric tensors which depend on
genus-0 and genus-1 data.

The universal equation corresponding to a formula due to Mumford is equivalent to
$\rho_{2;1} = 0$ (cf. \cite{G2}). Another equation due to Getzler \cite{G2} is equivalent to
$\rho_{2;2} = 0$. The equation $\rho_{2;3} = 0$ is due to Belorousski and Pandharipande
\cite{BP}. Therefore $\rho_{2;i} \in U$ for $i=1, 2, 3$.

\subsection{String and dilaton equations}
\label{sec:StringDilaton}

There are also universal equations which do not come from relations in the
tautological ring of the moduli space of stable curves. For example, the string equation,
dilaton equation and divisor equation. In this paper, we will only need the string and
dilaton equations.

The {\it string equation} has the form
\[ \gwiig{\vs} = \frac{1}{2} \delta_{g, 0} \eta_{\alpha \beta}
t_{0}^{\alpha} t_{0}^{\beta} \]
for $g \geq 0$. Since $ \nabla_{\vw} \vs = - \tau_{-}(\vw)$, taking covariant derivative
on both sides of the string equation, we obtain
\[
 \gwiig{ \vs \, \vw_{1} \, \cdots \, \vw_{k} }
    = \sum_{i=1}^{k} \gwiig{ \vw_{1} \, \cdots \,
            \left\{ \tau_{-}(\vw_{i}) \right\}
            \, \cdots \, \vw_{k} }
            + \delta_{g, 0} \nabla^{k}_{\vw_{1}, \cdots, \vw_{k}}
                \left( \frac{1}{2} \eta_{\alpha \beta}
            t_{0}^{\alpha} t_{0}^{\beta} \right).
\]
for any vector fields $\vw_{1}, \ldots, \vw_{k}$.
The {\it dilaton vector field} is defined by
\[ \vd = T({\vs}). \]
The {\it dilaton equation} has the form
\[ \left<\left< {\cal D} \right>\right>_{g} =
        (2g-2) F_{g} + \frac{1}{24} \, \chi(V)\delta_{g, 1} \]
for $g \geq 0$.
Since
$ \nabla_{\vw} \vd = - \vw$, taking covariant derivative
on both sides of the dilaton equation, we obtain
\[ \gwiig{\vd \, \vw_{1} \cdots \, \vw_{k}}
= (k+2g-2) \gwiig{\vw_{1} \, \cdots \vw_{k}}. \]

\section{Relations among lower genus universal equations}

\subsection{Relation among genus-2 universal equations}

Relations among genus-2 universal equations were studied in \cite{L2}.
Using the string and dilaton equations, it was proved that
\[ \rho_{2;2}(\vs, \vw) = \rho_{2;1}(\tau_{-}(\vw)) \]
for any vector field $\vw$ (cf. \cite[Theorem 2.6]{L2}).
Since $\tau_{-}$ is surjective, this shows that
the universal equation $\rho_{2;1} =0 $ follows from $\rho_{2;2} = 0$.
On the other hand,
it is straightforward to show that
\begin{eqnarray}
&& (\nabla_{T(\vw_{1})} \rho_{2;1})(\vw_{2}) - \rho_{2;2}(\vw_{1}, T(\vw_{2})) \nonumber \\
&=& \left\{ 3 \rho_{0}(\vw_{2}, \vw_{1}, \gua)  - \rho_{0}(\vw_{1}, \vw_{2}, \gua) \right\}
        \gwiitwo{T(\ga)}
    - \rho_{0}(\vw_{2}, T(\vw_{1}), \gua) \gwiitwo{\ga} \nonumber  \\
&&    + A_{2}(\vw_{1}, T(\vw_{2})) - (\nabla_{T(\vw_{1})} A_{1})(\vw_{2}) \label{eqn:MtoG}
    \end{eqnarray}
for any vector fields $\vw_{1}$ and $\vw_{2}$.
Two conclusions can be drawn from this equation: First, we obtain a genus-1 universal
relation
\begin{equation} \label{eqn:A1A2}
 A_{2}(\vw_{1}, T(\vw_{2})) - \left(\nabla_{T(\vw_{1})} \, A_{1} \right)(\vw_{2}) =0
\end{equation}
for any vector fields $\vw_{1}$ and $\vw_{2}$. Secondly, the universal equation
$\rho_{2;2}(\vw_{1}, T(\vw_{2})) = 0$ follows from $\rho_{2;1}=0$, the genus-0 topological
recursion relation and the genus-1 \eqref{eqn:A1A2}.

Note that a special case of the genus-0 topological recursion relation is the following
\[ T(\vw_{1}) \qp \vw_{2} = \rho_{0}(\vw_{1}, \vw_{2}, \gua) \ga =  0 \]
for all vector fields $\vw_{1}$ and $\vw_{2}$.
Taking covariant derivative of this equation with respect to $\vw_{3}$, we have
\[ \gwii{T(\vw_{1}) \, \vw_{2} \, \vw_{3} \, \gua} \ga  = \vw_{1} \qp \vw_{2} \qp \vw_{3} \]
for all vector fields $\vw_{1}$, $\vw_{2}$, and $\vw_{3}$.
Using these two equations, we obtain
\begin{eqnarray}
&& \rho_{2;3}(\vw_{1}, \vw_{2}, T(\vw_{3}))
    + (\nabla_{\vw_{1} \qp \vw_{2}} \rho_{2;1}) (\vw_{3})
    - \rho_{2;2} (\vw_{1} \qp \vw_{2}, \vw_{3}) \nonumber \\
&=& - B(\vw_{1}, \vw_{2}, T(\vw_{3})) - (\nabla_{\vw_{1} \qp \vw_{2}} A_{1}) (\vw_{3})
        + A_{2}(\vw_{1} \qp \vw_{2}, \vw_{3}). \label{eqn:MGtoBP}
\end{eqnarray}
Three conclusions can be drawn from this equation: First, we obtain a genus-1
universal equation
\begin{equation} \label{eqn:A1A2B}
 B(\vw_{1}, \vw_{2}, T(\vw_{3})) -
         A_{2}(\vw_{1} \qp \vw_{2}, \vw_{3})
         + (\nabla_{\vw_{1} \qp \vw_{2}} A_{1}) (\vw_{3}) = 0
\end{equation}
for all vector fields $\vw_{1}$, $\vw_{2}$, and $\vw_{3}$.
Secondly, the universal equation
\[\rho_{2;3}(\vw_{1}, \vw_{2}, T(\vw_{3}))=0 \]
follows from $\rho_{2;1} = 0$, $\rho_{2;2} = 0$, the genus-0 topological recursion relation,
and the genus-1  \eqref{eqn:A1A2B}. Thirdly, since any vector field $\vw$ can be written as
\[ \vw = \vs \qp \vw + T(\tau_{-}(\vw)), \]
combining \eqref{eqn:MGtoBP} and \eqref{eqn:MtoG}, we obtain the result that
the universal equation $\rho_{2;2}=0$ follows from $\rho_{2;3}=0$ and $\rho_{2:1} = 0$
together with genus-0 and genus-1 universal equations.
This is precisely \cite[Theorem 2.9]{L2}.

\subsection{Relation among genus-1 universal equations}

In this section, we will study the relations between genus-1 universal equations in
section \ref{sec:g1univ} and equations (\ref{eqn:A1A2}) and
(\ref{eqn:A1A2B}).
The main result is the following
\begin{thm} \label{thm:NoNewG1}
Equations (\ref{eqn:A1A2}) and
(\ref{eqn:A1A2B}) can be derived from known genus-0 and genus-1 relations
described in section \ref{sec:g0univ} and \ref{sec:g1univ}.
\end{thm}
 To prove this theorem,
we need first use genus-0 and genus-1 topological recursion relations and their
covariant derivatives to get rid of the $T$ operator in equations (\ref{eqn:A1A2}) and
(\ref{eqn:A1A2B}). The first three covariant derivatives of the genus-0 topological recursion
$\rho_{0} = 0$ gives
\begin{equation} \label{eqn:4ptTg0}
 \gwii{T(\vw_{1}) \, \vw_{2} \, \vw_{3} \, \vw_{4}}
    = \gwii{(\vw_{1} \qp \vw_{2}) \, \vw_{3} \, \vw_{4}},
\end{equation}
\begin{eqnarray}
 \gwii{T(\vw_{1}) \, \vw_{2} \, \ldots \, \vw_{5}}
& = & \gwii{ \{\vw_{1} \bullet \vw_{2} \} \, \vw_{3} \, \vw_{4} \, \vw_{5}}
\nonumber \\
&& + \gwii{ \{\vw_{1} \bullet \vw_{3} \} \, \vw_{2} \, \vw_{4} \, \vw_{5}}
\nonumber \\
&&  + \gwii{ \vw_{1} \, \vw_{2} \, \vw_{3} \, \{\vw_{4} \bullet \vw_{5} \}},
\label{eqn:5ptTg0}
\end{eqnarray}
and
\begin{eqnarray}
 \gwii{T(\vw_{1}) \, \vw_{2} \, \ldots \, \vw_{6}}
& = & \sum_{i=2}^{4}
\gwii{ \{\vw_{1} \bullet \vw_{i} \} \, \vw_{2} \, \ldots \, \widehat{\vw_{i}}
    \, \ldots \, \vw_{6}}
\nonumber \\
&&  + \gwii{ \vw_{1} \, \ldots \, \vw_{4} \,
         \{\vw_{5} \bullet \vw_{6}\} }
\nonumber \\
&& + \gwii{ \vw_{1} \, \vw_{2} \, \vw_{3} \, \gub}
    \gwii{ \gb \, \vw_{4} \, \vw_{5} \, \vw_{6} }
\nonumber \\
&& + \gwii{ \vw_{1} \, \vw_{2} \, \vw_{4} \, \gub}
    \gwii{ \gb \, \vw_{3} \, \vw_{5} \, \vw_{6} }
\nonumber \\
&& + \gwii{ \vw_{1} \, \vw_{3} \, \vw_{4} \, \gub}
    \gwii{ \gb \, \vw_{2} \, \vw_{5} \, \vw_{6} }.
\label{eqn:6ptTg0}
\end{eqnarray}
It's obvious that the left hand sides of these equations are symmetric with respect to
all arguments except $\vw_{1}$. This symmetry property is not apparent for the right hand
sides of these equations. Therefore when apply these equations to get rid of operator $T$,
we may get two seemingly different but actually equal expressions. For the convenience
of presentation, in this section we
will follow the rule that when applying these equations,
we always take the form which gives the smallest number of terms. A typical example
is the following: When applying \eqref{eqn:6ptTg0} to
$\gwii{ T(\vw) \, \vv \, \ga \, \gua \, \gb \, \gub}$, we choose $\vw_{1} = \vw$,
$\vw_{2} = \ga$, $\vw_{3} = \gua$, $\vw_{4} = \gb$, $\vw_{5} = \gub$, $\vw_{6} = \vv$.
Keeping in mind that we are summing over repeated indices, using apparent symmetries,
we obtain
\begin{eqnarray}
&& \gwii{ T(\vw) \, \vv \, \ga \, \gua \, \gb \, \gub} \nonumber \\
&=& 3 \gwii{ (\vw \qp \ga) \, \gua \, \gb \, \gub \, \vv}
    + \gwii{ \vw \, \ga \, \gua \, \gb \, (\gub \qp \vv)}  \nonumber \\
&&    + \gwii{\vw \, \ga \, \gua  \, \gum} \gwii{\gm \, \gb \, \gub \, \vv}
    + 2 \gwii{\vw \, \ga \, \gb \, \gum} \gwii{\gm \, \gua \, \gub \, \vv}.
    \label{eqn:typex1}
\end{eqnarray}
If we instead choose $\vw_{1} = \vw$,  $\vw_{2} = \vv$,
$\vw_{3} = \ga$, $\vw_{4} = \gua$, $\vw_{5} = \gb$, $\vw_{6} = \gub$, we obtain
\begin{eqnarray}
&& \gwii{ T(\vw) \, \vv \, \ga \, \gua \, \gb \, \gub} \nonumber \\
&=& \gwii{ (\vw \qp \vv) \, \ga \, \gua \, \gb \, \gub}
    + 2 \gwii{ (\vw \qp \ga) \, \gua \, \gb \, \gub \, \vv}
    + \gwii{ \vw \, \vv \, \ga \, \gua \, (\gb \qp \gub)}  \nonumber \\
&&    + 2 \gwii{\vw \, \vv \, \gua  \, \gum} \gwii{\gm \, \ga \, \gb \, \gub}
    +  \gwii{\vw \, \ga \, \gua \, \gum} \gwii{\gm \, \vv \, \gb \, \gub}.
    \label{eqn:typex2}
\end{eqnarray}
Although the right hand sides of \eqref{eqn:typex2} and \eqref{eqn:typex1}
look very different from each other, they are actually equal since the left hand sides
of these equations are the same. The right hand side of \eqref{eqn:typex2} has 5 terms.
It is slightly more complicated than the right hand side of \eqref{eqn:typex1} which has only
4 terms. So following the rule stated earlier, we should use \eqref{eqn:typex1}
rather than \eqref{eqn:typex2}.

To eliminate operator $T$ in genus-1 correlation functions, we need the first
three covariant derivatives of the genus-1 topological recursion relation
$\rho_{1}=0$, which have the following forms:
\begin{equation}
\gwiione{T({\cal W}) \, \vv} \, \, = \, \,
    \gwiione{ \{\vw \bullet \vv\} } +
     \frac{1}{24} \gwii{{\cal W} \, \vv \, \gum  \, \gm},
\label{eqn:derTRRg1}
\end{equation}
\begin{eqnarray}
\gwiione{T({\cal W}) \, \vv_{1} \, \vv_{2}}
& = &   \gwiione{ \{\vw \bullet \vv_{1}\} \, \vv_{2} } +
    \gwiione{ \{\vw \bullet \vv_{2}\} \, \vv_{1} } +
    \gwii{ \vw \, \vv_{1} \, \vv_{2} \, \gua } \gwiione{\ga}
\nonumber \\
&&  + \frac{1}{24} \gwii{{\cal W} \, \vv_{1} \, \vv_{2} \, \gum  \, \gm},
\label{eqn:2derTRRg1}
\end{eqnarray}
and
\begin{eqnarray}
&& \gwiione{T({\cal W}) \, \vv_{1} \, \vv_{2} \, \vv_{3}} \nonumber \\
& = &   \sum_{i=1}^{3} \left\{
        \gwiione{ (\vw \qp \vv_{i}) \, \vv_{1} \, \ldots \,
                        \widehat{\vv_{i}} \, \ldots  \vv_{3} }
                        + \gwii{ \vw \, \vv_{1} \, \ldots \,
                        \widehat{\vv_{i}} \, \ldots  \vv_{3} \,
                                \gua} \gwiione{\ga \, \vv_{i}} \right\}
        \nonumber \\
&& +  \gwii{ \vw \, \vv_{1} \, \vv_{2} \, \vv_{3} \, \gua } \gwiione{\ga}
  + \frac{1}{24} \gwii{{\cal W} \, \vv_{1} \, \vv_{2} \, \vv_{3} \, \gua  \, \ga}
\label{eqn:3derTRRg1}
\end{eqnarray}
for all vector fields. Note that both sides of these equations are symmetric with respect
to vector fields $\vv_{i}$'s. There is no ambiguity on how to apply these formulas.

Besides the genus-0 and genus-1 topological recursion relations and their covariant
derivatives, we also need universal equations $C_{0} = 0$ and $G=0$
(see Section \ref{sec:g0univ} and \ref{sec:g1univ} for the definitions of tensors $C_{0}$
and $G$). The application of $C_{0} = 0$ is straightforward. Most times we only need
a slightly weaker version of this equation, i.e. the associativity of the quantum product.
However it is less obvious how to apply the covariant derivatives of $C_{0}=0$.
This turns out to be one of the most subtle points in the proof of Theorem~\ref{thm:NoNewG1}.
In the proof, we will include tensors $\nabla^{i} C_{0}$ in the relations among genus-1
universal tensors. It is then apparent where and how covariant derivatives of $C_{0}=0$
are used.

The tensors for the covariant derivatives of $C_{0}$ are defined by
\[ C_{i}(\vw_{1}, \cdots, \vw_{i+4}) :=
    (\nabla_{\vw_{i+4}} C_{i-1})(\vw_{1}, \cdots, \vw_{i+3}) \]
for $i \geq 1$. These tensors lie in $U$, the ideal of tensors given universal equations.
We need the first three covariant derivatives of $C_{0}$ in the proof of
Theorem~\ref{thm:NoNewG1}. Using \eqref{eqn:derQP},  we can compute these tensors
more explicitly. For example,
\begin{eqnarray*}
 C_{1}(\vw_{1}, \cdots, \vw_{5})
&=& \gwii{ (\vw_{1} \qp \vw_{2}) \, \vw_{3} \, \vw_{4} \, \vw_{5}}
    + \gwii{ \vw_{1} \, \vw_{2} \, (\vw_{3} \qp \vw_{4}) \, \vw_{5}} \\
&&  - \gwii{(\vw_{1} \qp \vw_{3}) \, \vw_{2} \, \vw_{4} \, \vw_{5}}
        - \gwii{ \vw_{1} \, \vw_{3} \, (\vw_{2} \qp \vw_{4}) \, \vw_{5}}
\end{eqnarray*}
and
\begin{eqnarray*}
&& C_{2}(\vw_{1}, \cdots, \vw_{6})  \\
&=& \gwii{ (\vw_{1} \qp \vw_{2}) \, \vw_{3} \, \cdots \, \vw_{6}}
    + \gwii{ \vw_{1} \, \vw_{2} \, (\vw_{3} \qp \vw_{4})
                \, \vw_{5} \, \vw_{6}} \\
&&  - \gwii{(\vw_{1} \qp \vw_{3}) \, \vw_{2} \, \vw_{4}
                \, \vw_{5} \, \vw_{6}}
        - \gwii{ \vw_{1} \, \vw_{3} \, (\vw_{2} \qp \vw_{4})
                        \, \vw_{5} \, \vw_{6}} \\
&& + \gwii{\vw_{1} \, \vw_{2} \, \vw_{5} \, \gua}
      \gwii{\ga \, \vw_{3} \, \vw_{4} \, \vw_{6}}
 + \gwii{\vw_{1} \, \vw_{2} \, \vw_{6} \, \gua}
      \gwii{\ga \, \vw_{3} \, \vw_{4} \, \vw_{5}} \\
&& - \gwii{\vw_{1} \, \vw_{3} \, \vw_{5} \, \gua}
      \gwii{\ga \, \vw_{2} \, \vw_{4} \, \vw_{6}}
 - \gwii{\vw_{1} \, \vw_{3} \, \vw_{6} \, \gua}
      \gwii{\ga \, \vw_{2} \, \vw_{4} \, \vw_{5}}
\end{eqnarray*}
for any vector fields $\vw_{1}, \ldots, \vw_{6}$. The third derivative $C_{3}$
also has a similar expression, but more complicated. So we omit the explicit
expression for $C_{3}$.

We are now ready to prove Theorem~\ref{thm:NoNewG1}.

{\bf Proof of Theorem~\ref{thm:NoNewG1}}:
First using genus-0 and genus-1 topological recursion relations and their derivatives,
we can get rid of operator $T$ in the expression
\[ A_{2}(\vw_{1}, T(\vw_{2}))
    - \left(\nabla_{T(\vw_{1})} \, A_{1} \right)(\vw_{2}). \]
We then check that the resulting
expression is equal to certain combinations of tensors $G$ and $C_{1}$.
More precisely, modulo
\[ I_{a}[\rho_{0}, \, \nabla \rho_{0}, \, \nabla^{2} \rho_{0}, \, \nabla^{3} \rho_{0}, \,\,
     \rho_{1}, \, \nabla \rho_{1}, \, \nabla^{2} \rho_{1}, \,
      C_{0}], \]
we have
\begin{eqnarray}
&& 120 \left\{ A_{2}(\vw_{1}, T(\vw_{2}))
    - \left(\nabla_{T(\vw_{1})} \, A_{1} \right)(\vw_{2}) \right\}  \nonumber \\
&=& - \frac{1}{24} G(\vw_{1}, \vw_{2}, \gua, \ga)
    + \left\{\frac{1}{6} C_{1}(\vw_{2}, \ga, \vw_{1}, \gua, \gb)
            \right. \nonumber \\
&& \left. + \frac{19}{2} C_{1}(\vw_{1}, \gb, \ga, \gua, \vw_{2})
    - \frac{21}{2} C_{1}(\vw_{2}, \gb, \gua, \ga, \vw_{1})
    \right\} \gwiione{\gub}
    \label{eqn:A1A2G}
\end{eqnarray}
for any vector fields $\vw_{1}$ and $\vw_{2}$.
This shows that \eqref{eqn:A1A2} follows from genus-0 and genus-1
universal equations in Section~\ref{sec:univ}.
The proof of \eqref{eqn:A1A2G} is straightforward, but somehow tedious.

\Eqref{eqn:A1A2B} can be treated in a similar fashion. We first use genus-0 and
genus-1 topological recursion relations and their derivatives to get rid of operator
$T$ in the expression
\[ B(\vw_{1}, \vw_{2}, T(\vv)) -
         A_{2}(\vw_{1} \qp \vw_{2}, \vv)
         + (\nabla_{\vw_{1} \qp \vw_{2}} A_{1}) (\vv), \]
then check that the resulting
expression is equal to certain combinations of tensors $G$, $\nabla G$, $C_{1}$,
$C_{2}$, and $C_{3}$.
More precisely, modulo
\[ I_{a}[\rho_{0}, \, \nabla \rho_{0}, \, \nabla^{2} \rho_{0}, \, \nabla^{3} \rho_{0}, \,\,
     \rho_{1}, \, \nabla \rho_{1}, \, \nabla^{2} \rho_{1}, \, \nabla^{3} \rho_{1}, \,
      C_{0}], \]
 we have
\begin{eqnarray}
&& 720 \, \left\{ B(\vw_{1}, \vw_{2}, T(\vv)) -
         A_{2}(\vw_{1} \qp \vw_{2}, \vv)
         + (\nabla_{\vw_{1} \qp \vw_{2}} A_{1}) (\vv) \right\} \nonumber \\
&=& 24 \, G(\vv, \vw_{1}, \vw_{2}, \gua) \gwiione{\ga}
    +  (\nabla_{\ga} G)(\vw_{1}, \vw_{2}, \vv, \gua) \nonumber \\
&&  +\frac{1}{4} \, \left\{ (\nabla_{\vv} G)(\vw_{1}, \vw_{2}, \ga, \gua)
    -(\nabla_{\vw_{1}} G)(\vv, \vw_{2}, \ga, \gua)
    - (\nabla_{\vw_{2}} G)(\vv, \vw_{1}, \ga, \gua) \right\}
    \nonumber     \\
&&   - 3 \, \gwiione{\vv \, \gb} \,
    \left\{ C_{1}(\vw_{1}, \vw_{2}, \ga, \gua, \gub)
            -C_{1}(\vw_{1}, \ga, \vw_{2}, \gub, \gua) \right. \nonumber \\
&& \hspace{200pt} \left.  -C_{1}(\vw_{2}, \ga, \vw_{1}, \gub, \gua) \right\} \nonumber \\
&& + \gwiione{\vw_{1} \, \gb}
    \left\{ 27 \, C_{1}(\vw_{2}, \ga, \vv, \gub, \gua)
            - 34 \,  C_{1}(\vv, \ga, \vw_{2}, \gub, \gua) \right. \nonumber \\
&& \hspace{200pt} \left.  - 7 \, C_{1}(\vw_{2}, \gub, \ga, \gua, \vv) \right\} \nonumber \\
&& + \gwiione{\vw_{2} \, \gb}
    \left\{ 27 \, C_{1}(\vw_{1}, \ga, \vv, \gub, \gua)
            - 34 \,  C_{1}(\vv, \ga, \vw_{1}, \gub, \gua) \right. \nonumber \\
&& \hspace{150pt} \left.  - 7 \, C_{1}(\vw_{1}, \gub, \ga, \gua, \vv) \right\}
        \nonumber \\
&& - 8 \, \gwiione{\ga \, \gb}
    \left\{ 2 \, C_{1}(\vw_{1}, \vw_{2}, \gua, \gub, \vv)
            -  C_{1}(\vw_{1}, \vv, \vw_{2}, \gua, \gub) \right. \nonumber \\
&& \hspace{200pt} \left.  - C_{1}(\vw_{2}, \vv, \vw_{1}, \gua, \gub) \right\}
            \nonumber \\
&& - 120 \, \gwiione{\ga} \gwiione{\gb}
    \left\{ 3 \, C_{1}(\vw_{1}, \vw_{2}, \gua, \gub, \vv)
            -  C_{1}(\vw_{1}, \vv, \vw_{2}, \gua, \gub) \right. \nonumber \\
&& \hspace{200pt} \left.  - C_{1}(\vw_{2}, \vv, \vw_{1}, \gua, \gub) \right\} \nonumber  \\
&& - \gwiione{\gb} \left\{
    19 \, C_{2}(\vw_{2}, \gub, \vw_{1}, \ga, \vv, \gua)
    + 19 \, C_{2}(\vw_{1}, \gub, \vw_{2}, \ga, \vv, \gua) \right. \nonumber \\
&& \hspace{60pt}
    + 19 \, C_{2}(\vw_{1}, \vw_{2}, \ga, \gua, \vv, \gub)
    - 12 \, C_{2}(\vv, \ga, \gub, \gua, \vw_{1}, \vw_{2}) \nonumber \\
&& \hspace{60pt}
    - 3 \, C_{2}(\vw_{1}, \vv, \vw_{2}, \gub, \ga, \gua)
    - 3 \, C_{2}(\vw_{2}, \vv, \vw_{1}, \gub, \ga, \gua) \nonumber \\
&& \hspace{60pt} \left.
    - 28 \, C_{2}(\vw_{1}, \ga, \vw_{2}, \gub, \vv, \gua)
    - 28 \, C_{2}(\vw_{2}, \ga, \vw_{1}, \gub, \vv, \gua) \right\} \nonumber  \\
&& - C_{3}(\vw_{1}, \vw_{2}, \ga, \gua, \vv, \gb, \gub) \label{eqn:A1A2BG}
\end{eqnarray}
for any vector fields $\vw_{1}$, $\vw_{2}$, and $\vv$.
This implies that \eqref{eqn:A1A2B} follows from genus-0 and genus-1
universal equations in Section~\ref{sec:univ}.
The proof of \eqref{eqn:A1A2BG}  is very tedious but nevertheless straightforward.
All what is needed is to plug the expressions for the relevant tensors
into both sides of \eqref{eqn:A1A2BG} and check that they are equal after
obvious cancellations.
$\Box$

%%%%%%%%%%%%%%%%%%%%%%%%%%%%%%%%%%%%%%%%%%%%%%%%%%%%%

\vspace{30pt} \noindent
Department of Mathematics  \\
University of Notre Dame \\
Notre Dame,  IN  46556, USA \\

\vspace{10pt} \noindent E-mail address: {\it xliu3@nd.edu}


\begin{thebibliography}{399}

\bibitem[BP]{BP} Belorousski, P. and Pandharipande, R.,
        {\it A descendent relation in genus 2},
    Ann. Scuola Norm. Sup. Pisa Cl. Sci. (4) 29 (2000) 171-191.
\bibitem[DZ]{DZ} Dubrovin, B., Zhang, Y.,
           {\it Bihamiltonian hierarchies in 2D topological field
                theory at one-loop approximation},
           Comm. Math. Phys. 198 (1998) 311 - 361.
\bibitem[EHX]{EHX} Eguchi, T., Hori, K., and Xiong, C.,
        {\it Quantum Cohomology and Virasoro Algebra},
        Phys. Lett. B402 (1997) 71-80.
\bibitem[G1]{G1} Getzler, E.,
        {\it Intersection theory on $\bar{M}_{1,4}$ and elliptic
                Gromov-Witten Invariants},
        J. Amer. Math. Soc. 10 (1997) 973-998
\bibitem[G2]{G2} Getzler, E.,
        {\it Topological recursion relations in genus 2},
       Integrable systems and algebraic geometry (Kobe/kyoto, 1997)
    73-106.
\bibitem[KM]{KM} Kontsevich, M. and Manin, Y.,
        {\it Relations between the correlators of the topological sigma-model
        coupled to gravity}, Comm. Math. Phys., 196 (1998) 385-398.
\bibitem[LiT]{LiT} Li, J. and Tian, G.,
        {\it Virtual moduli cycles and Gromov-Witten invariants of
                general symplectic manifolds},
        Topics in symplectic 4-manifolds (Irvine, CA, 1996), 47-83.
\bibitem[L1]{L1} Liu, X.,
    {\it Elliptic Gromov-Witten invariants and Virasoro conjecture},
    Comm. Math. Phys. 216 (2001), 705-728.
\bibitem[L2]{L2} Liu, X.,
    {\it Quantum product on the big phase space and Virasoro conjecture},
    to appear in Advances in Mathematics ({\textsf math.AG/0104030}).
\bibitem[LT]{LT} Liu, X. and Tian, G.,
        {\it Virasoro constraints for quantum cohomology},  \\
        J. Diff. Geom. 50 (1998), 537 - 591.
\bibitem[M]{M} Manin, Y.,
    {\it Frobenius manifolds, quantum cohomology, and moduli spaces},
    Amer. Math. Soc. Colloq. Publ. 47, Providence, RI, 1999.
\bibitem[RT]{RT} Ruan, Y. and  Tian, G.,
        {\it Higher genus symplectic invariants and sigma models coupled
         with gravity}, Invent. Math. 130 (1997), 455-516.
\bibitem[W]{W} Witten, E.,
        {\it Two dimensional gravity and intersection theory on
                Moduli space},
        Surveys in Diff. Geom., 1 (1991), 243-310.


\end{thebibliography}
\end{document}